\documentclass[11pt, a4paper]{article}
\usepackage[russian, english]{babel}
\usepackage{amsmath,amsfonts,amsthm,amssymb,amscd}
\numberwithin{equation}{section}
\date{}
\newtheorem{teo}{Theorem}[section]
\newtheorem{lem}[teo]{Lemma}

\newtheorem{rem}[teo]{Remark}
\newtheorem{prop}[teo]{Proposition}
\begin{document}
\title{On the Two Obstacles Problem in Orlicz-Sobolev Spaces and Applications}
\author{\textit{Dedicated to Professor V.V. Zhikov on the occasion of his 70th anniversary}\\ \ \\
\textrm{Jos$\acute{e}$ Francisco Rodrigues\footnote{rodrigue@fc.ul.pt}\quad\textrm{and } Rafayel Teymurazyan\footnote{rafayel@ptmat.fc.ul.pt}}\\ \ \\
\textit{Centro de Matem$\acute{a}$tica e Aplica\c{c}\~{o}es Fundamentais (CMAF/FCUL)}\\
\textit{Universidade de Lisboa,}\\
\textit{Av. Prof. Gama Pinto, 2}\\
\textit{1649-003 Lisboa, Portugal}
}
\maketitle
\begin{abstract}We prove the Lewy-Stampacchia inequalities for the two obstacles problem in
abstract form for T-monotone operators. As a consequence for a general class of quasi-linear
elliptic operators of Ladyzhenskaya-Uraltseva type, including p(x)-Laplacian type
operators, we derive new results of $C^{1,\alpha}$ regularity for the solution. We also
apply those inequalities to obtain new results to the N-membranes problem and the regularity and monotonicity properties to obtain the existence of a solution to a quasi-variational
problem in (generalized) Orlicz-Sobolev spaces.
\end{abstract}
\section{Introduction}
\quad

We consider the two obstacles problem for monotone operators (possibly degenerate or singular) of the type

\begin{equation}
Au=-\textrm{div}(a(x,|\nabla u|)\nabla u)
\end{equation}
with a Dirichlet boundary condition in a bounded domain $\Omega\subset\mathbb{R}^n$.\\

The two obstacles problem for the operator $(1.1)$ consists of finding  $u\in\mathbb{K}_\psi^\varphi$ such that

\begin{equation}
\int_\Omega a(x,|\nabla u|)\nabla u\cdot\nabla(v-u)dx\geq\int_\Omega f(v-u)dx,\textrm{ }\forall v\in\mathbb{K}_\psi^\varphi,
\end{equation}
where
\begin{equation}
\mathbb{K}_\psi^\varphi=\{v\in W_0^{1,G}(\Omega):\textrm{ }\psi\leq
v\leq\varphi\textrm{ a.e. in }\Omega,\}
\end{equation}
$\psi,\varphi\in W^{1,G}(\Omega)$ (for the definition of the Orlicz-Sobolev space $W^{1,G}(\Omega)$ see next section), where $G$ is related to $a$ by

\begin{equation}
G(x,t)=\int_0^t a(x,s)sds,\textrm{ }x\in\Omega,\textrm{ }t\geq0.
\end{equation}

$W_0^{1,G}(\Omega)$ is the closure of $C_0^\infty(\Omega)$ in $W^{1,G}(\Omega)$ and is a Banach space of weakly differentiable functions $v$ with

$$
\int_\Omega G(x,|\nabla v(x)|)dx<\infty.
$$

In particular, when in $(1.1)$ we have $a(x,t)=t^{p(x)-2}$,
with $p(x)>1$ a given bounded function in $\Omega$, we deal
with problems involving variable growth conditions, the so called $p(x)$-Laplacians. The study of such problems has been stimulated by problems in elasticity
(see [32]), in fluid dynamics (see [3], [9], [29], [33]), image processing models [8] and problems in the calculus of variations with $p(x)$-growth conditions
(see [2], [22], [23], [32], [34]) and some more general class of differential operators
(see [4], [6], [12], [33]).\\

Here we are specially interested in the more general class of quasi-linear operators of Ladyzhenskaya-Uraltseva type (see [19], [20]), when $a(x,t):\Omega\times\mathbb{R}^+\rightarrow\mathbb{R}$ is given by a function measurable and bounded in $x$ for all $t>0$ and Lipschitz continuous in $t$, a.e. $x\in\Omega$, and, such that, there are positive constants $\underline{a}<\overline{a}$
\begin{equation}
0<\underline{a}\leq\frac{ta_t(x,t)}{a(x,t)}+1\leq\overline{a}\textrm{ for }t>0,
\end{equation}
where $a_t=\partial a/\partial t$, and also $\lim_{t\rightarrow0+}ta(x,t)=0$.

The assumption $(1.5)$, in fact, implies
\begin{equation}
\big(a(x,|\xi|)\xi-a(x,|\zeta|)\zeta\big)\cdot(\xi-\zeta)>0,\textrm{ }\forall \xi,\zeta\in\mathbb{R}^n,\textrm{ }\xi\neq \zeta
\end{equation}
and $\textrm{ and }\lim_{t\rightarrow\infty}ta(x,t)=\infty$ for a.e. $x\in\Omega$ (see [7], for instance).

As a consequence, we have the uniqueness of the solution and also the weak maximum principle for $A$. In this work, after recalling the natural functional framework of the Orlicz-Sobolev spaces associated with $A$, in section 2, we extend some classical properties on the bilateral problem $(1.2)-(1.3)$ in this more general framework, for instance, by including examples like $a(x,t)=\alpha(x)t^{p(x)-2}\log(\beta(x)t+\gamma(x))$ with bounded functions $\gamma(x)$, $p(x)>1$, and $\alpha(x)$, $\beta(x)>0$ a.e. in $x\in\Omega$.

In section 3, we use the continuity property of the truncation operator $v\mapsto v^+=v\vee0=\sup(v,0)$ for the strong topology of $W^{1,G}(\Omega)$ to extend some continuous dependence results in $W_0^{1,G}(\Omega)$ of the variational solutions to $(1.2)-(1.3)$ with respect to the data.

In section 4, we prove the Lewy-Stampacchia inequalities (here $a\wedge b=\inf(a,b)$ and $a\vee b=\sup(a,b)$)

\begin{equation}
A\varphi\wedge f\leq Au\leq A\psi\vee f
\end{equation}
in an abstract form, extending the approach of Mosco [25] to the two obstacles problem, that includes the above class of operators. Although the inequalities $(1.7)$ are known, in particular, for linear operators (see [30]), our proof is new and more general. As a consequence, under additional (H$\ddot{o}$lder) continuity hypothesis on $x\mapsto a(x,\cdot)$, we obtain the same regularity for the solution $u$ of the two obstacles problem as in the equation without constraints (see [20]). For instance, for bounded obstacles $\varphi$, $\psi$ we conclude that $u\in C^{1,\alpha}(\Omega)$, if we impose $f$, $A\varphi$, $A\psi\in L^\infty(\Omega)$, a regularity obtained in [19] with different assumptions.

Finally, in section 5, we give two new applications to systems of obstacle type. In the case of the N-membranes problem, when $\textbf{u}=(u_1,\ldots,u_N)$ has the constraint

$$
u_1\geq u_2\geq\ldots\geq u_N\textrm{ a.e. in }\Omega,
$$
we extend some of the results of [5], in particular, the $C^{1,\alpha}$ regularity and the strong approximation in $(W_0^{1,G}(\Omega))^N$ by solutions of a penalized system. For the case of a special class of implicit double obstacle problems, when the obstacles depend on the solution in the form

$$
\bigvee_{i\neq j}(u_j-\psi_{ij})\leq u_i\leq\bigwedge_{i\neq j}(u_j+\varphi_{ij}), \textrm{ }i=1,\ldots,N,
$$
where $\varphi_{ij}$, $\psi_{ij}$ are certain given positive constants, we are able to show the existence of a minimal and maximal solution for the corresponding system of quasi-variational inequalities, which is of the type arising in problems of stochastic impulse control (see [14], [25], [30] and [31]).

\section{Preliminaries on Orlicz-Sobolev Spaces}

The Orlicz spaces and the Orlicz-Sobolev spaces as well are defined for the $\textit{Young functions}$ (or $\textit{N-functions}$) like $G$ defined by $(1.4)$. Let for all $x\in\Omega$, $g(x,\cdot):\mathbb{R}\rightarrow\mathbb{R}$ is an odd, increasing homeomorphism from $\mathbb{R}$ onto $\mathbb{R}$; $g(x,t)>0$, when $t>0$, while the function $G:{\Omega}\times\mathbb{R}\rightarrow\mathbb{R}$,
$$
G(x,t):=\int_{0}^{t}g(x,s)ds
$$
for all $x\in{\Omega}$ and all $t\geq0$ belongs to class $G$ (see [26], p. 33), i.e. $G$ satisfies the following conditions:

$(i)$ for all $x\in\Omega$, $G(x,\cdot):[0,\infty)\rightarrow\mathbb{R}$ is an increasing function, $\lim_{t\rightarrow\infty}g(x,t)=\infty$, $G(x,0)=0$ and $G(x,t)>0$ whenever $t>0$;

$(ii)$ for every $t\geq0$, $G(\cdot,t):\Omega\rightarrow\mathbb{R}$ is a measurable function.\\

For a Young function $G$, we define the $\textit{(generalized) Orlicz class}$,
$$
K_G(\Omega)=\{u:\Omega\rightarrow\mathbb{R},\textrm{ measurable;
}\rho_G(u):=\int_\Omega G(x,|u(x)|)dx<\infty\}
$$
and also the $\textit{(generalized) Orlicz space}$,
$$
L^G(\Omega)=\{u:\Omega\rightarrow\mathbb{R}, \textrm{measurable;} \lim_{\lambda\rightarrow0^+}\int_\Omega G(x,\lambda|u(x)|)dx=0\}.
$$
which is a Banach space endowed with the $\textit{Luxemburg norm}$
$$
\|u\|_G=\inf\bigg\{\mu>0; \int_\Omega G\bigg(x, \frac{|u(x)|}{\mu}\bigg)dx\leq1\bigg\}.
$$
The $\textit{(generalized) Orlicz-Sobolev space}$ is defined as follows
$$
W^{1,G}(\Omega)=\bigg\{u\in L^G(\Omega);\frac{\partial u}{\partial
x_i}\in L^G(\Omega), i=1,2,...,n\bigg\},
$$
is also a Banach space with the norm:
$$
\|u\|_{1,G}=\||\nabla u|\|_G+\|u\|_G.
$$
These spaces are more general that the usual Lebesgue or Sobolev spaces, but many properties of functions in these spaces can be extended. In particular, the Poincar$\acute{e}$ type inequality
$$
\int_\Omega G(x,|u|)dx\leq\int_\Omega G(x,c|\nabla u|)dx,
$$
holds for any $u\in W_0^{1,G}(\Omega)$, where $c$ is twice the diameter of
$\Omega$ (see [15]).\\
The H\"{o}lder inequality extends to (see [26], Theorem 13.13)
$$
\bigg|\int_\Omega uvdx\bigg|\leq
C\|u\|_G\|v\|_{\overline{G}},\textrm{ }\forall u\in
L^G(\Omega), v\in L^{\overline{G}}(\Omega),
$$
where $C$ is a positive constant, and $\overline{G}$ is the $\textit{conjugate}$ Young function of $G$, that is,
$$
\overline{G}(x,t)=\sup_{s>0}\{ts-G(x,s); s\in\mathbb{R}\},\textit{ }\forall x\in\Omega,\textrm{ }t\geq0.
$$
We also denote by $G^*$ the $\textit{Sobolev conjugate}$ of $G$, that is
$$
(G^*)^{-1}(x,t)=\int_0^t\frac{G^{-1}(x,s)}{s^{(n+1)/n}}ds,
$$
provided
$$
\int_1^\infty\frac{G^{-1}(x,s)}{s^{(n+1)/n}}ds=\infty.
$$
It is also well known, that
$$
L^\infty(\Omega)\hookrightarrow L^G(\Omega)\hookrightarrow L^1(\Omega),
$$
with continuous imbeddings, and if $\Omega$ is a bounded domain with a smooth boundary, then the imbedding $W_0^{1,G}(\Omega)\hookrightarrow L^{G^*}(\Omega)$ is continuous.

In this work given the function
$a:\Omega\times\mathbb{R}^+\rightarrow\mathbb{R}^+$ satisfying the assumption $(1.5)$. Let the mapping $g:\Omega\times\mathbb{R}\rightarrow\mathbb{R}$ be defined by
\begin{equation}
g(x,t):=\left\{\begin{array}{ll}
a(x,|t|)t,\textrm{ if  }t\neq0,\\
0, \textrm{ if }t=0.
\end{array}\right.
\end{equation}
Then $g$ satisfies the conditions $(i)-(ii)$, and the corresponding $G$ is a strictly convex Young function.

In our case, $(1.5)$ implies
\begin{equation}
0<1+\underline{a}\leq\frac{tg(x,t)}{G(x,t)}\leq\overline{a}+1\textrm{ ,a.e. } x\in\Omega\textrm{ }t\geq0,
\end{equation}
and $G$ satisfies the so called $\Delta_2$-$\textit{condition}$ (see, for instance [1]), which implies that $L^G(\Omega)=K_G(\Omega)$ (Theorem 8.13 in [26]), and that $L^G(\Omega)$ is a linear separable space.

Relation $(2.2)$ assures (see Proposition 2.2 in [23]) that $L^G(\Omega)$ is an uniformly convex space and thus,
a reflexive space.\\

In Orlicz-Sobolev spaces an important role is played by $\rho_G(u)$ - the $\textit{modular}$ of the $L^G(\Omega)$ space. If $u_m$, $u\in L^G(\Omega)$ then (see [23] for the proofs) when $\|u\|_G>1$ it holds $\|u\|_G^{1+\underline{a}}\leq\rho_G(u)\leq\|u\|_G^{1+\overline{a}}$, and we have
$$
\|u_m\|_G\rightarrow\infty\textrm{ }\Leftrightarrow\textrm{ }\rho_G(u_m)\rightarrow\infty,
$$
and when $\|u\|_G<1$ also $\|u\|_G^{1+\overline{a}}\leq\rho_G(u)\leq\|u\|_G^{1+\underline{a}}$, which implies
$$
\|u_m-u\|_G\rightarrow0\textrm{ }\Leftrightarrow\textrm{ }\rho_G(u_m-u)\rightarrow0.
$$

We refer to [10], [11], [18], [26] for further properties of (generalized) Lebesgue-Sobolev spaces.

Here we need the following lemma:

\begin{lem}
If $u\in W^{1,G}(\Omega)$, then $u^+,u^-\in W^{1,G}(\Omega)$ and
$$
Du^+=\left\{\begin{array}{ll}
Du, &\textrm{ if }u>0,\\
0, &\textrm{ if }u\leq0,
\end{array}\right.
\text{ and }
Du^-=\left\{\begin{array}{ll}
0, &\textrm{ if }u\geq0,\\
-Du, &\textrm{ if }u<0.
\end{array}\right.
$$
\end{lem}
Here $u^+=u\vee0$, $u^-=-u\wedge0$. This lemma holds in $W_0^{1,G}(\Omega)$ as well, and therefore these Orlicz-Sobolev spaces are closed with respect to
$$
u\vee v=u+(v-u)^+\textrm{ and }u\wedge v=u-(u-v)^+.
$$
The proof of this lemma is due to Gossez (see [16]) and is basically the same as for usual Sobolev spaces (see, for instance, Theorem 1.56 in [30]).
\begin{lem} The embeddings
$W_0^{1,\overline{a}+1}(\Omega)\hookrightarrow W_0^{1,G}(\Omega)\hookrightarrow W_0^{1,\underline{a}+1}(\Omega)$ are continuous.
\end{lem}
\textbf{Proof.} The second part was observed in [20]. In order to prove the first part, by Theorem 8.12 in [1] it is enough to check, that $t^{\overline{a}+1}$ dominates $G$ near infinity, which is true, since by taking the log in the right inequality of $(2.2)$, we conclude that there exists $c>0$ and $T>0$ such that $G(x,t)\leq ct^{\overline{a}+1}$ for $t>T$.$\blacksquare$
\begin{prop}
If $u_m\rightarrow u$ in $W^{1,G}(\Omega)$, then $u_m^+\rightarrow u^+$ in $W^{1,G}(\Omega)$.
\end{prop}
\textbf{Proof.} This follows by convexity of $G$ and by the inequality (see [22]):
$$
G(x,t+s)\leq 2^{1+\underline{a}}(2+\underline{a})\big(G(x,t)+G(x,s)\big)\textrm{ a.e. }a\in\Omega\textrm{ }\forall t,s>0.
$$
We have
$$
\int_\Omega G(x,|\nabla(u_m^+-u^+)|)dx=\int_\Omega G(x,|\chi_{\{u_m>0\}}\nabla u_m-\chi_{\{u>0\}}\nabla u|)dx
$$
$$
\leq c\bigg(\int_\Omega\chi_{\{u_m>0\}} G(x,|\nabla(u_m-u)|)dx+\int_\Omega|\chi_{\{u_m>0\}}-\chi_{\{u>0\}}|G(x,|\nabla u|)dx\bigg)\rightarrow0,
$$
and so $\nabla u_m^+\rightarrow\nabla u^+$ in $L^G(\Omega)$. Arguing in the same way, we get also $u_m^+\rightarrow u^+$ in $L^G(\Omega)$, which completes the proof.$\blacksquare$
\begin{rem}
a) Assuming $G(x,t)=G(t)$, i.e. $G$ is independent of variable $x$, we say that $L^G$
and $W^{1,G}$ are Orlicz spaces, respectively Orlicz-Sobolev spaces (see [1]).\\
b) Assuming $G(x,t)=|t|^{p(x)}$ with $p(x)\in L^\infty(\Omega)$, $p(x)\geq\underline{p}>1$ a.e. in $\Omega$, we denote $L^G(\Omega)$ by $L^{p(x)}(\Omega)$ and $W^{1,G}(\Omega)$ by $W^{1,p(x)}(\Omega)$ and we refer them as variable exponents Lebesgue spaces, respectively
variable exponents Sobolev spaces.\\
c) Our framework enables us to work with spaces which are more general than those described in a) and b). Besides the example given in the introduction for $a(x,t)=\alpha(x)t^{p(x)-2}\log(\beta(x)t+\gamma(x))$ with $p(x)$, $\gamma(x)>1$ and $\alpha(x)$, $\beta(x)>0$ a.e. $x\in\Omega$, we could consider any linear combination with positive coefficients or any composition of functions satisfying a condition like $(1.5)$.
\end{rem}

\section{Variational solutions}

We introduce the energy functional $J:W_0^{1,G}(\Omega)\rightarrow\mathbb{R}$ by
\begin{equation}
J(u)=\int_\Omega G(x,|\nabla u|)\textrm{ }dx,\textrm{ }\forall u\in W_0^{1,G}(\Omega)
\end{equation}
which is strictly convex, weakly lower semi-continuous and coercive in $W_0^{1,G}(\Omega)$ (see [21]). Moreover, $J$ is G$\hat{a}$teaux differentiable, and $J'(u)=Au$ at $u\in W_0^{1,G}(\Omega)$ is given by (see [23], for instance)
\begin{equation}
\langle Au,v\rangle=\int_\Omega a(x,|\nabla u|)\nabla u\cdot\nabla v\textrm{ }dx,\textrm{ }v\in W_0^{1,G}(\Omega).
\end{equation}
Hence $Au\in W^{-1,G}(\Omega)=\big(W_0^{1,G}(\Omega)\big)'$, the topological dual of $W_0^{1,G}(\Omega)$ and if we assume $f\in L^{\overline{G}^*}(\Omega)\subset W^{-1,G}(\Omega)$, $\overline{G}^*$ being the conjugate Young function of the Sobolev conjugate of $G$, we can rewrite the problem $(1.3)-(1.4)$ in the form:\\
find $u\in\mathbb{K}_\psi^\varphi$, such that,
\begin{equation}
\langle Au-L,v-u\rangle\geq0\textrm{ }\forall v\in\mathbb{K}_\psi^\varphi,
\end{equation}
where we set
$$
\langle L,v\rangle=\int_\Omega fv\textrm{ }dx,\textrm{ }\forall v\in W_0^{1,G}(\Omega).
$$
\begin{prop}
Under the condition $(1.6)$ the operator $A$ is strictly T-monotone, i.e., for any $u,v\in W^{1,G}(\Omega)$
$$
\langle Au-Av,(u-v)^+\rangle>0,\textrm{ if }0\neq(u-v)^+\in W_0^{1,G}(\Omega).
$$
\end{prop}
\textbf{Proof.} In fact, using Lemma 2.1, by $(1.6)$ we have
$$
\int_\Omega\bigg\{\big[a(x,|\nabla u|)\nabla u-a(x,|\nabla v|)\nabla v\big]\cdot\nabla(u-v)^+\bigg\}dx
$$
$$
=\int_{\{u>v\}}\bigg\{\big[a(x,|\nabla u|)\nabla u-a(x,|\nabla v|)\nabla v\big]\cdot(\nabla u-\nabla v)\bigg\}dx>0,
$$
if $(u-v)^+\neq0$, i.e. if $meas\{u>v\}>0$.$\blacksquare$\\

In this section we assume $(1.5)$ and
\begin{equation}
\varphi,\textrm{ }\psi\in W^{1,G}(\Omega)\textrm{ such that }\mathbb{K}_\psi^\varphi\neq\emptyset,
\end{equation}
for which it is sufficient to assume $\varphi\geq\psi$ a.e. in $\Omega$ and both $\varphi^-$, $\psi^+\in W_0^{1,G}(\Omega)$.

\begin{teo}
The problem $(3.3)$ has a unique solution $u=u(f,\varphi,\psi)\in\mathbb{K}_\psi^\varphi$ and is equivalent to minimize in $\mathbb{K}_\psi^\varphi$ the functional $F$, defined by
\begin{equation}
F(v)=\int_\Omega G(x,|\nabla v|)\textrm{ }dx-\int_\Omega fv\textrm{ }dx,\textrm{ }v\in W_0^{1,G}(\Omega).
\end{equation}
Moreover, if $\hat{u}$ denotes the solution corresponding to $\hat{f}$, $\hat{\varphi}$, $\hat{\psi}$, then
$$
f\geq\hat{f},\textrm{ }\varphi\geq\hat{\varphi},\textrm{ }\psi\geq\hat{\psi}\textrm{ implies }u\geq\hat{u}\textrm{ a.e. in }\Omega.
$$
\end{teo}

\textbf{Proof.} The existence and uniqueness are standard results for strictly monotone, coercive and potential operators, as observed in more general Orlicz-Sobolev spaces already in [17]. The comparison property follows easily by the T-monotonicity (see, for instance, [25] or [27]): take $v=u\vee\hat{u}\in\mathbb{K}_\psi^\varphi$ in $(3.3)$ and $v=u\wedge\hat{u}\in\mathbb{K}_{\hat{\psi}}^{\hat{\varphi}}$ in $\hat{(3.3)}$ for $\hat{u}$. By addition, one finds
$$
\langle A\hat{u}-Au,(\hat{u}-u)^+\rangle+\langle L-\hat{L},(\hat{u}-u)^+\rangle\leq0.
$$
Since $L-\hat{L}\geq0$, and $A$ is strictly T-monotone, one immediately deduces $(\hat{u}-u)^+=0$, which means that $u\geq\hat{u}$.$\blacksquare$
\begin{rem}
This argument also shows a weak maximum principle in the form: if $Au\geq A\hat{u}$ in $\Omega$ and $u\geq\hat{u}$ on $\partial\Omega$ in the sense $(\hat{u}-u)^+\in W_0^{1,G}(\Omega)$, then $u\geq\hat{u}$ in $\Omega$.
\end{rem}
Similarly we have a ``$L^\infty$-continuous dependence'' property, even without the $L^\infty$ regularity on the solutions:
\begin{prop}
$$
\|u-\hat{u}\|_{L^\infty(\Omega)}\leq\|\varphi-\hat{\varphi}\|_{L^\infty(\Omega)}
\vee\|\psi-\hat{\psi}\|_{L^\infty(\Omega)},
$$
where $u$ and $\hat{u}$ are the corresponding solutions of the problem $(3.3)$ and $\hat{(3.3)}$ with the same $f$.
\end{prop}
\textbf{Proof.} Let $l=\|\varphi-\hat{\varphi}\|_{L^\infty(\Omega)}
\vee\|\psi-\hat{\psi}\|_{L^\infty(\Omega)}<\infty$. Set $v=u+(\hat{u}-u-l)^+\in\mathbb{K}_\psi^\varphi$ in $(3.3)$ and $\hat{v}=\hat{u}-(\hat{u}-u-l)^+\in\mathbb{K}_{\hat{\psi}}^{\hat{\varphi}}$ in $\hat{(3.3)}$.
By addition one gets
$$
I_1:=\int_\Omega\big(a(x,|\nabla\hat{u}|)\nabla\hat{u}-a(x,|\nabla u|)\nabla u\big)\cdot\nabla\big(\hat{u}-u-l\big)^+\leq0.
$$
On the other hand, recalling $(1.6)$ we know, that if $meas\{\hat{u}>u+l\}>0$
$$
I_2:=\int_{\{\hat{u}>u+l\}}\big(a(x,|\nabla\hat{u}|)\nabla\hat{u}-a(x,|\nabla u|)\nabla u\big)\cdot(\nabla\hat{u}-\nabla u)^+>0.
$$
Since $I_1=I_2$, we conclude $\hat{u}-u-l\leq0$. Reversing the role of $u$ with $\hat{u}$ we get $u-\hat{u}-l\leq0$.$\blacksquare$\\

Exactly as in Proposition 4.5 of [23], we have the following interesting result.
\begin{lem}
Let $(1.5)$ holds and $u_m\rightharpoonup u$ weakly in $W_0^{1,G}(\Omega)$. Then if \begin{equation}
\limsup_{m\rightarrow\infty}\langle Au_m, u_m-u\rangle\leq0
\end{equation}
then $u_m\rightarrow u$ in $W_0^{1,G}(\Omega)$ strongly, and $Au_m\rightarrow Au$ in $W^{-1,G}(\Omega)$ strongly.
\end{lem}
\begin{teo}
Under the assumption $(1.5)$ let $u_m$ be the solution to $(3.3)$ with compatible data $(f_m,\varphi_m,\psi_m)$, such that
$$
f_m\rightarrow f\in L^{\overline{G}^*},\textrm{ }\varphi_m\rightarrow\varphi\textrm{ and }\psi_m\rightarrow\psi\textrm{ in }W^{1,G}(\Omega)\textrm{ strongly}.
$$
Then $u_m\rightarrow u$ strongly in $W_0^{1,G}(\Omega)$, where $u$ is the solution of the limit problem $(3.3)$ with $(f,\varphi,\psi)$.
\end{teo}
\textbf{Proof.} For arbitrary $v\in\mathbb{K}_\psi^\varphi$ we obtain that $v_m=(v\wedge\varphi_m)\vee\psi_m\in\mathbb{K}_{\psi_m}^{\varphi_m}$ and, by Proposition 2.3, $v_m\rightarrow v$ in $W_0^{1,G}(\Omega)$. Since, by coerciveness, $\|u_m\|_{1,G}\leq C$, there is $u\in W_0^{1,G}(\Omega)$, such that, for a subsequence $u_m\rightharpoonup u$ in $W_0^{1,G}(\Omega)$ weakly and in $L^1(\Omega)$. So $u\in\mathbb{K}_\psi^\varphi$ and by lower semi-continuity we have
$$
F(u)\leq\liminf_{m\rightarrow\infty}F(u_m)\leq\lim_{m\rightarrow\infty}F(v_m)=F(v),\textrm{ }\forall v\in\mathbb{K}_\psi^\varphi.
$$
Hence $u$ solves $(3.3)$ and then, since $w_m=(u\wedge\varphi_m)\vee\psi_m\rightarrow u$ in $W_0^{1,G}(\Omega)$,
$$
\limsup_{m\rightarrow\infty}\langle Au_m,u_m-u\rangle\leq\lim_{m\rightarrow\infty}\langle Au_m, u_m-w_m\rangle+\lim_{m\rightarrow\infty}\langle Au_m,w_m-u\rangle
$$
$$
\leq\lim_{m\rightarrow\infty}\int_\Omega f_m(w_m-u_m)=0.
$$
By Lemma 3.5, this implies the strong convergence of the whole sequence $u_m$ to $u$.$\blacksquare$
\begin{rem}
In the case of only one obstacle problem, lower $(\psi)$ or upper $(\varphi)$ obstacle (corresponding to take formally $\varphi=+\infty$ or $\psi=-\infty$ respectively) all previous results of this section hold in similar way. For instance, Proposition 3.4 takes the form
$$
\|u-\hat{u}\|_{L^\infty(\Omega)}\leq\|\psi-\hat{\psi}\|_{L^\infty(\Omega)}\textrm{ or }\|u-\hat{u}\|_{L^\infty(\Omega)}\leq\|\varphi-\hat{\varphi}\|_{L^\infty(\Omega)}
$$
respectively.
\end{rem}

\section{Lewy-Stampacchia inequalities and its consequences}
\quad
Let $X$ be a real reflexive Banach space, which is a lattice with respect to a partial order $"\leq"$, and $V$ is a sublattice of $X$ (i.e. $V$ contains the $\sup$ and $\inf$ of its elements). Then for every $v\in X$ one can write $v=v^+-v^-$, where $v^+=v\vee0$, $v^-=-v\wedge0$ are the positive and negative parts of
$v$ respectively. Thus, $X$ is generated by the cone $P$, $P-P=X$,
$$
P=\{v\in X:\textrm{ }v\geq0\}.
$$
The subspace of $V'$ (the dual space of $V$) generated by the cone
$$
P'=\{v'\in V':\textrm{ }\langle v',v\rangle\geq0,\textrm{ }\forall
v\in P\},
$$
is called the $\textbf{dual order}$ of $V$ and denoted by $V^*$, i.e. $V^*=P'-P'$.
We shall consider an operator $A:X\rightarrow V'$ with the properties:\\
a) $\textbf{hemi-continuous}$, i.e. the mapping $t\rightarrow\langle A(u+tv),w\rangle$ is continuous on $[0,1]$, $u,v\in X$, $w\in V$;\\
b) $\textbf{coercive}$, i.e. $\exists w_0\in X$ such that
$$
\textrm{if }u-w_0\in V,\textrm{ }\lim_{\|u\|\rightarrow\infty}\frac{\langle Au,u-w_0\rangle}{\|u\|}=\infty,\textrm{ }u\in X,
$$
c) $\textbf{strictly T-monotone}$, i.e.
$$
\langle Au-Av, (u-v)^+\rangle>0
$$
for all $u,v\in X$ such that $0\neq(u-v)^+\in V$.\\

We recall that $L\geq M$ for $L$, $M\in V'$, if $L-M\in P'$, i.e. $L-M$ is positive on the positive elements of $V$.  Then, if $A$ is strictly T-monotone, it satisfies a weak comparison principle: if $Au\leq Av$, then $u\leq v$.\\

In this framework we recall the following interesting fact (see [25]).
\begin{teo}
Let $A:X\rightarrow V'$ be a strictly T-monotone, coercive and hemi-continuous operator. Let also $u,v\in V$ be such that $Au,Av\in V^*$. Then $A(u\wedge v)$, $A(u\vee v)\in V^*$ and $A(u\wedge v)\geq Au\wedge Av$, $A(u\vee v)\leq Au\vee Av$ for the dual order in $V'$.\\
Consequently, if $Au_i\in V^*$, then for $\forall i\in I\subset\{1,2,\ldots,N\}$
$$
A\bigg(\bigwedge_{i\in I}u_i\bigg)\geq\bigwedge_{i\in I}Au_i\textrm{ and }A\bigg(\bigvee_{i\in I}u_i\bigg)\leq\bigvee_{i\in I}Au_i.
$$
\end{teo}
In fact, we can take $X=W^{1,G}(\Omega)$, $V=W_0^{1,G}(\Omega)$ and the operator $A$ defined by $(3.2)$ under our assumptions.\\

Now we prove the Lewy-Stampacchia inequalities in abstract form.
\begin{teo}
Let $X$ be a real reflexive ordered Banach space, $V$ be a closed subspace of $X$, which is a sublattice of $X$, $A:X\rightarrow V'$ be a strictly T-monotone, coercive and hemi-continuous operator. Let also there are given two elements $\varphi$,$\psi\in X$, $\psi\leq\varphi$, and $L\in V'$. If\\
(a)
$$
\exists\textrm{ }\Lambda\in V'\textit{ such that } \Lambda\geq
L\textit{ and }\Lambda\geq A\psi\textit{ in } V'
$$
$$
\exists\textrm{ }\lambda\in V'\textit{ such that } \lambda\leq
L\textit{ and }\lambda\leq A\varphi\textit{ in } V',
$$
(b)
$$
(\psi-v)^{+}\in V\textrm{ and }(v-\varphi)^{+}\in V,\textit{ }\forall v\in V,
$$
and let $u$ be the solution of
\begin{equation}
u\in\mathbb{K}_{\psi}^{\varphi}:\textit{ }\langle Au-L,v-u\rangle\geq0,\textit{ }\forall
v\in\mathbb{K}_{\psi}^{\varphi},
\end{equation}
where
\begin{equation}
\mathbb{K}_{\psi}^{\varphi}=\{u\in V: \psi\leq u\leq\varphi\}.
\end{equation}
Then the following dual estimates hold
\begin{equation}
\lambda\leq Au\leq\Lambda\textrm{ in } V'.
\end{equation}
In particular, if $L, A\varphi, A\psi\in V^*$, then also $Au\in
V^*$, and $(4.3)$ gives
\begin{equation}
L\wedge A\varphi\leq Au\leq L\vee A\psi\textrm{ in } V^{*}.
\end{equation}
\end{teo}
\textbf{Proof.} The assertion $(4.4)$ assertion holds, since, if $L, A\varphi,A\psi\in V^{*}$, one can take $\lambda=L\wedge A\varphi=L-(L-A\varphi)^{+}$ and $\Lambda=L\vee A\psi=L+(A\psi-L)^{+}$ in $(4.3)$.\\
To prove the upper bound of $(4.3)$, consider the unique solution $z\in V$ of auxiliary variational inequality
\begin{equation}
z\leq u:\textit{ }\langle Az-\Lambda,w-z\rangle\geq0,\textit{ }\forall w\in V,\textit{ } w\leq u.
\end{equation}
It is enough to show, that $z=u$, since then taking $w=u-v$ in $(4.5)$ for an arbitrary $v\geq0$, it readily follows
\begin{displaymath}
Au-\Lambda=Az-\Lambda\leq0\textrm{ in } V'.
\end{displaymath}
To prove that $z=u$, let us first prove that $z\geq\psi$. Recalling the first condition of $(b)$ and taking $w=z+(\psi-z)^{+}=\psi\vee z\leq u$ in $(4.5)$, we get
\begin{displaymath}
\langle\Lambda-Az,(\psi-z)^{+}\rangle\leq0.
\end{displaymath}
Hence, since $\Lambda\geq A\psi$, one obtains
\begin{displaymath}
\langle A\psi-\Lambda,(\psi-z)^{+}\rangle+\langle\Lambda-Az,(\psi-z)^{+}\rangle\leq0
\end{displaymath}
which, by the strict T-monotonicity of $A$, implies, that $(\psi-z)^{+}=0$. This means that $z\geq\psi$.\\
Let us now prove that $z\geq u$. Since $u$ solves $(4.1)$ and in $(4.5)$ $z\leq u$, one has that $z\leq\varphi$. In other words, $z\in\mathbb{K}_{\psi}^{\varphi}$. Take $w=z\vee u=z+(u-z)^{+}\leq u$ in $(4.5)$ and $v=u\wedge z=u-(u-z)^{+}\geq\psi$ in $(4.1)$. Then, by addition and recalling the fact $\Lambda\geq L$ from the conditions $(a)$, we conclude
\begin{displaymath}
\langle Au-Az,(u-z)^{+}\rangle\leq\langle L-\Lambda,(u-z)^{+}\rangle\leq0,
\end{displaymath}
and, since $A$ is strictly T-monotone, it follows that $(u-z)^{+}=0$ and so $z\geq u$.\\
So, the unique solution $u$ of $(4.1)$ is also the unique solution of $(4.5)$. We already know, that this implies the upper bound in $(4.3)$.

To prove the lower bound in $(4.3)$, consider the unique solution $z\in V$ of auxiliary variational inequality
\begin{equation}
z\geq u:\textit{ }\langle Az-\lambda,w-z\rangle\geq0,\textit{ }\forall w\in V,\textit{ } w\geq u.
\end{equation}
The steps are similar as above. It is enough to show, that $z=u$, since then taking $w=u+v$ in $(4.6)$ for an arbitrary $v\geq0$, it readily follows
\begin{displaymath}
\lambda-Au=\lambda-Az\leq0\textrm{ in } V'.
\end{displaymath}
To prove that $z=u$, let us first prove that $z\leq\varphi$. Recalling the second condition in $(b)$ and taking $w=z-(z-\varphi)^{+}=z\wedge\varphi\geq u$ in $(4.6)$, we get
\begin{displaymath}
\langle Az-\lambda,(z-\varphi)^{+}\rangle\leq0.
\end{displaymath}
Hence, since $\lambda\leq A\varphi$, one obtains
\begin{displaymath}
\langle\lambda-A\varphi,(z-\varphi)^{+}\rangle+\langle Az-\lambda,(z-\varphi)^{+}\rangle\leq0,
\end{displaymath}
which, by the strict T-monotonicity of $A$, implies, that $(z-\varphi)^{+}$=0, i.e. $z\leq\varphi$.\\
Let us now prove that $z\leq u$. Since $u$ solves $(4.1)$ and in $(4.6)$ $z\geq u$, one has that $z\geq\psi$. In other words, $z\in\mathbb{K}_{\psi}^{\varphi}$. Take $w=z\wedge u=z-(z-u)^{+}\geq u$ in $(4.6)$ and $v=u\vee z=u+(z-u)^{+}\leq\varphi$ in $(4.1)$. Then, by subtracting and recalling that $\lambda\leq L$ from $(a)$, we conclude
\begin{displaymath}
\langle Az-Au,(z-u)^{+}\rangle\leq\langle\lambda-L,(z-u)^{+}\rangle\leq0,
\end{displaymath}
and, since $A$ is strictly T-monotone, it follows, that $(z-u)^{+}=0$ and so $z\leq u$.

So, the unique solution $u$ of $(4.1)$ is also the unique solution of $(4.5)$ (we already know, that this implies the upper bound in $(4.3)$) and the unique solution of $(4.6)$, which, as we see, implies the lower bound in $(4.3)$. $\blacksquare$

\begin{rem} Theorem 4.2 is still true when $X$ is just an ordered Banach space and the operator $A$ is strictly T-monotone, provided the problem $(4.1)$, $(4.2)$ has solution. This result extends Theorem 4.1 of [25] that was restricted to the abstract one obstacle problem, and Theorem 4.35 of [30] was stated only for linear second order operators in Sobolev spaces, which was recently extended to the $p$-Laplacian in [28] and to general Leray-Lions operators in [24].
\end{rem}

Let $u$ be the unique solution of $(4.1)$, $\underline{u}$ be the unique
solution of $(4.1)$ in $\mathbb{K}_{\psi}=\{v\in V: \psi\leq
v\}$ (i.e. $\underline{u}$ is the solution of one (lower) obstacle
problem) and let $\overline{u}$ be the unique solution of $(4.1)$ in
$\mathbb{K}^{\varphi}=\{v\in V: v\leq\varphi\}$ (i.e. $\overline{u}$ is
the solution of one (upper) obstacle problem). We can recover the Lewy-Stampacchia inequalities for the one obstacle problem easily in the following way.
\begin{prop}
i) If $A\varphi\geq L$, then $u=\underline{u}$, and one can rewrite $(4.4)$ as
$$
L\leq Au\leq L\vee A\psi\textrm{ in } V^{*}.
$$
ii) If $A\psi\leq L$, then $u=\overline{u}$, and one can rewrite $(4.4)$ as
$$
L\wedge A\varphi\leq Au\leq L\textrm{ in } V^{*}.
$$
\end{prop}

\textbf{Proof.} $i)$ By taking $\underline{v}=\underline{u}-(\underline{u}-\varphi)^+=\underline{u}\wedge\varphi\in\mathbb{K}_{\psi}$ in the lower-obstacle problem
$$
\underline{u}\in\mathbb{K}_{\psi}\textrm{ : }\langle A\underline{u}-L,\underline{v}-\underline{u}\rangle\geq0,\textrm{ }\forall\underline{v}\in\mathbb{K}_{\psi},
$$
we get
$$
\langle A\underline{u}-L,(\underline{u}-\varphi)^+\rangle\leq0.
$$
On the other hand, by assumption,
$$
\langle L-A\varphi,(\underline{u}-\varphi)^+\rangle\leq0.
$$
Recalling the strictly T-monotonicity of $A$, by addition, we conclude $(\underline{u}-\varphi)^+=0$ and so $\underline{u}\leq\varphi$, which means, that $\underline{u}\in\mathbb{K}_{\psi}^{\varphi}$. The uniqueness of the solution gives $\underline{u}=u$.

$ii)$ Similarly, by taking $\overline{v}=\overline{u}+(\psi-\overline{u})=\overline{u}\vee\psi\in\mathbb{K}^{\varphi}$ in the upper-obstacle problem
$$
\overline{u}\in\mathbb{K}_{\varphi}\textrm{ : }\langle A\overline{u}-L,\overline{v}-\overline{u}\rangle\geq0,\textrm{ }\forall\overline{v}\in\mathbb{K}_{\varphi},
$$
we get, using
$$
\langle A\psi-A\overline{u},(\psi-\overline{u})^+\rangle=\langle A\psi-L,(\psi-\overline{u})^+\rangle+\langle L-A\overline{u},(\psi-\overline{u})^+\rangle\leq0.
$$
So $(\psi-\overline{u})^+=0$, which means that $\overline{u}\in\mathbb{K}_{\psi}^{\varphi}$, and we conclude by the uniqueness of the solution that $\overline{u}=u$.$\blacksquare$\\

We have already observed that the nonlinear operator $A$ given by $(1.1)$ under the generalized condition $(1.5)$ of Ladyzhenskaya and Uraltseva type satisfies the hypothesis of Theorem 3.2 in $W_0^{1,G}(\Omega)$. In particular, if we also assume the existence of $C>0$, such that
\begin{equation}
|f(x)|\leq C,\textrm{ }A\psi\leq C\textrm{ and }A\varphi\geq-C\textrm{ a.e. }x\in\Omega,
\end{equation}
the solution of $(3.3)$ satisfies
$$
-C\leq(f\wedge A\varphi)(x)\leq Au(x)\leq(f\vee A\psi)(x)\leq C\textrm{ a.e. }x\in\Omega,
$$
and we conclude that
$$
Au\in L^\infty(\Omega).
$$
Hence, the regularity of the solution of the two obstacles problem (and, similarly, for each one obstacle problem as well) is the same as bounded solutions of the respective equation without obstacles, as it was established in [20]. So if, we assume, in addition to the continuous differentiability in $t$ of $a(x,t)$, that for some constants $C_\beta>0$ and $0<\beta\leq1$
\begin{equation}
|a(x,t)t-a(y,t)t|\leq C_\beta(1+a(x,t)t)|x-y|^\beta,\textrm{ }\forall x,y\in\overline{\Omega},\textrm{ }t>0,
\end{equation}
from [19], [20] we immediately conclude the following interesting regularity result (which for the $p(x)$-Laplacian operator was re-established in [13]).
\begin{teo}
Under the assumptions $(1.5)$, $(3.4)$, $(4.7)$ and $(4.8)$ the solution $u$ of the obstacle problem $(3.3)$ with $\varphi$, $\psi\in L^\infty(\Omega)$ is $C^{1,\alpha}(\Omega)$, for some $0<\alpha<1$. If also $\partial\Omega\in C^{1,\alpha}$, then $u\in C^{1,\alpha}(\overline{\Omega})$.
\end{teo}
\begin{rem}
The H$\ddot{o}$lder continuity of the gradient to the one or two obstacles problem has been obtained by Lieberman [19] when the obstacles are also $C^{1,\alpha}$. However, the condition $(4.7)$ on $\varphi$ and $\psi$ do not imply that $\varphi$, $\psi$ belong also to $C^{1,\alpha}$. So the result of Theorem 4.5 does not follow nor imply Lieberman's $C^{1,\alpha}$ regularity results.
\end{rem}

\section{Applications to Systems}

\subsection{N-membranes problem}
\quad For the operator $A$ the N-membranes problem consists of: to find $\textbf{u}=(u_1,u_2,...,u_N)\in\mathbb{K}_N$ satisfying
\begin{equation}
\sum_{i=1}^N\int_\Omega a(x,|\nabla
u_i|)\nabla u_i\cdot\nabla(v_i-u_i)dx\geq\sum_{i=1}^N\int_\Omega
f_i(v_i-u_i)dx,\textrm{ }\forall(v_1,...v_N)\in\mathbb{K}_N,
\end{equation}
where $\mathbb{K}_N$ is the convex subset of the Orlicz-Sobolev space
$[W_0^{1,G}(\Omega)]^N$, defined by
$$
\mathbb{K}_N=\{(v_1,...,v_N)\in[W_0^{1,G}(\Omega)]^N:\textrm{
}v_1\geq...\geq v_N\textrm{ a.e. in }\Omega\},
$$
and $f_1$,...,$f_N\in L^{\overline{G}^*}(\Omega)$. As in Theorem 3.2 the existence and uniqueness of the solution to $(5.1)$ follows easily.

\begin{teo}
If $a(x,t)$ satisfies to the conditions above (namely $(1.6)$), then the solution $(u_1,...,u_N)$ of the N-membranes problem for $A$ satisfies the
following Lewy-Stampacchia type estimates:
$$
f_1\wedge Au_1\leq f_1\vee...\vee f_N
$$
$$
f_1\wedge f_2\leq Au_2\leq f_2\vee...\vee f_N
$$
$$
\vdots
$$
$$
f_1\wedge...\wedge f_{N-1}\leq Au_{N-1}\leq f_{N-1}\vee f_N
$$
$$
f_1\wedge...\wedge f_N\leq Au_N\leq f_N
$$
a.e. in $\Omega$.
\end{teo}

\textbf{Proof.} Observe that choosing
$(v,u_2,...u_N)\in\mathbb{K}_N$, with $v\in\mathbb{K}_{u_2}$, we see
that $u_1\in\mathbb{K}_{u_2}$ solves the ``lower-obstacle problem'' with $f=f_1$, and so (recall Proposition 4.4)
$$
f_1\leq Au_1\leq f_1\vee Au_2\textrm{ a.e. in }\Omega.
$$
Analogously, we see that $u_j\in\mathbb{K}_{u_{j+1}}^{u_{j-1}}$
solves the two-obstacles problem with $f=f_j$,
$j=2,3,...N-1$, and satisfies, by $(4.4)$,
$$
f_j\wedge Au_{j-1}\leq Au_j\leq f_j\vee Au_{j+1}\textrm{ a.e. in
}\Omega.
$$
Since $u_N\in\mathbb{K}^{u_{N-1}}$, then (recall Proposition 4.4)
$$
f_N\wedge Au_{N-1}\leq Au_N\leq f_N\textrm{ a.e. in }\Omega.
$$
The proof concludes by simple iteration.$\blacksquare$\\

Since for each $i=1,2,\ldots,N$, the respective $Au_i$ maybe controlled in $L^\infty(\Omega)$ by
$$
\bigwedge_{1\leq j\leq i}f_j\leq Au_i\leq\bigvee_{i\leq j\leq N}f_j,
$$
if $a(x,t)$ is H$\ddot{o}$lder continuous in $x$ and continuously differentiable in $t$, by the regularity of [20] we also have
\begin{teo}
If $f_i\in L^\infty(\Omega)$ $(i=1,2,\ldots N)$ and $(1.5)$, $(4.8)$ hold,
then the solution $\textbf{u}$ of the N-membranes problem has
$C^{1,\alpha}(\Omega)$ regularity. If also
$\partial\Omega\in C^{1,\alpha}$, then the solution $\textbf{u}$ to N-membranes
problem belongs to $[C^{1,\alpha}(\overline{\Omega})]^N$ for some
$0<\alpha<1$.
\end{teo}

As in [5] we may also approximate the variational inequality using bounded penalization. Defining
$$
\xi_0=\max\bigg\{\frac{f_1+\cdots+f_i}{i}:i=1,\ldots,N\bigg\},
$$
$$
\xi_i=i\xi_0-(f_1+\cdots+f_i)\textrm{ for }i=1,\ldots,N,
$$
we observe that
$$
\left\{\begin{array}{ll}
\xi_i\geq0, &\textrm{ if }i\geq1,\\
(\xi_{i-1}-\xi_{i-2})-(\xi_i-\xi_{i-1})=f_i-f_{i-1}, & \textrm{ if }i\geq2.
\end{array}\right.
$$
For $\varepsilon>0$, let $\theta_\varepsilon$ be define as follows:
$$
\theta_\varepsilon:\overline{\mathbb{R}}\rightarrow\overline{\mathbb{R}},\textrm{ }s\mapsto\left\{\begin{array}{ll}
0, &\textrm{ if }s\geq0,\\
s/\varepsilon, & \textrm{ if }-\varepsilon<s<0,\\
-1, & \textrm{ if }s\leq-\varepsilon.
\end{array}\right.
$$
The approximate problem is given by the system for $u_i^\varepsilon\in W^{1,G}_0(\Omega)$
\begin{equation}
Au_i^\varepsilon+\xi_i\theta_\varepsilon(u_i^\varepsilon-u_{i+1}^\varepsilon)-
\xi_{i-1}\theta_\varepsilon(u_{i-1}^\varepsilon-u_i^\varepsilon)=f_i, \textrm{ in }\Omega,\textrm{ }i=1,\ldots,N,
\end{equation}
with the convention $u_0^\varepsilon=+\infty$, $u_{N+1}^\varepsilon=-\infty$.
Since the operator $A$ is strictly T-monotone, then arguing as in [5], we get:
\begin{teo}
If the operator $A$ satisfies the assumption $(1.5)$, then\\
$(i)$ the problem $(5.2)$ has a unique solution $(u_1^\varepsilon,\ldots,u_N^\varepsilon)\in [W_0^{1,G}(\Omega)]^N$. This solution satisfies
$$
u_i^\varepsilon\leq u_{i-1}^\varepsilon+\varepsilon\textrm{ for }i=2,\ldots,N.
$$
$(ii)$ $(u_1^\varepsilon,\ldots,u_N^\varepsilon)\rightarrow(u_1,\ldots,u_N)$ in $[W_0^{1,G}(\Omega)]^N$ strongly as $\varepsilon\rightarrow0$, where $(u_1,\ldots,u_N)$ is the solution of the N-membranes problem. Under assumptions of Theorem 5.2 the convergence holds in $[C^1(\overline{\Omega})]^N$.
\end{teo}

\textbf{Proof.} As observed in [5], the bounded penalization $B_\varepsilon$
$$
\langle B_\varepsilon\textbf{v},\textbf{w}\rangle=
\sum_{i=1}^{N}\int_\Omega\big[\xi_i\theta_\varepsilon(v_i-v_{i+1})-\xi_{i-1}\theta_\varepsilon(v_{i-1}-v_i)\big]w_i\textrm{ }dx
$$
is monotone in $[L^1(\Omega)]^N$. Therefore $(i)$ follows exactly in the same way as in Proposition 2.1 of [5].\\
Since $\langle B_\varepsilon\textbf{v},\textbf{v}-\textbf{u}^\varepsilon\rangle=0$ for any $\textbf{v}\in\mathbb{K}_N$, the weak convergence of $(ii)$ follows also exactly as in [5] by monotonicity arguments, since we have
\begin{equation}
\langle A\textbf{u}^\varepsilon,\textbf{v}-\textbf{u}^\varepsilon\rangle
\geq\langle\textbf{L},\textbf{v}-\textbf{u}^\varepsilon\rangle,\textrm{ }\forall\textbf{v}\in\mathbb{K}_N.
\end{equation}
Here we have used the notations
\begin{equation}
\langle A\textbf{w},\textbf{v}\rangle=\sum_{i=1}^N\langle Aw_i,v_i\rangle\textrm{ and }\langle\textbf{L},\textbf{v}\rangle=\sum_{i=1}^N\int_\Omega f_iv_i\textrm{ }dx.
\end{equation}
Since $\textbf{u}^\varepsilon\rightharpoonup\textbf{u}\in\mathbb{K}_N$ weakly in $[W_0^{1,G}(\Omega)]^N$, setting $(5.2)$ in variational form, using $(5.3)$ with $\textbf{v}=\textbf{u}$ and taking the $\limsup$ we obtain the strong convergence (using Lemma 3.5). Finally, if also $\textbf{f}\in[L^\infty(\Omega)]^N$ then the penalization term in $(5.2)$ is also bounded in $L^\infty$ and the $u_i^\varepsilon$ are uniformly bounded in $C^{1,\alpha}$.$\blacksquare$

\subsection{A Quasi-Variational Problem}

Another interesting application of the Lewy-Stampacchia inequalities
is its application in studying quasi-variational inequalities. Some of
these problems are related to a stochastic switching game.

In this section we assume that $\partial\Omega\in C^{1,\beta}$ for some $0<\beta<1$, and the operator $A$ has the $C^{1,\alpha}$ regularity property.

We consider the following quasi-variational problem (with the notations $(5.4)$):\\
find $\textbf{u}\in\mathbb{K}(\textbf{u})$, such that,
\begin{equation}
\langle A\textbf{u}-\textbf{L}, \textbf{v}-\textbf{u}\rangle\geq 0,\textrm{ }\forall\textbf{v}\in\mathbb{K}(\textbf{u}),
\end{equation}
where
$$
\mathbb{K}(\textbf{u})=\{\textbf{v}\in [W^{1,G}_0(\Omega)]^N:\Psi_i(\textbf{u})\leq v_i\leq\Phi_i(\textbf{u}),\textrm{ }i=1,2,\ldots,N\},
$$
$\textbf{f}=(f_1,\ldots,f_N)\in [L^{\overline{G}^*}(\Omega)]^N$, and for $\textbf{v}=(v_1,\ldots,v_N)$ we set
$$
\Phi_i(\textbf{v}):=\bigwedge_{i\neq j}(v_j+\varphi_{ij})\textrm{ }(\varphi_{ij}\textrm{ are positive constants}),
$$
$$
\Psi_i(\textbf{v}):=\bigvee_{i\neq j}(v_j-\psi_{ij})\textrm{ }(\psi_{ij}\textrm{ are positive constants}),
$$
for $i,j=1,\ldots,N$, i.e. we consider the problem $(5.1)$ with
$\mathbb{K}_N=\mathbb{K}(\textbf{u})$ (the obstacles themselves
depend on the solution). A similar problem for linear
operators (and for the one obstacle) was considered in [14] and [30].

We set
\begin{equation}
\mu:=\bigwedge_{i=1}^Nf_i,\textrm{ and }\nu:=\bigvee_{i=1}^Nf_i
\end{equation}
and we consider the unique solution $\underline{u}^0$ and $\overline{u}^0$ in $W_0^{1,G}(\Omega)$ solving, respectively,
\begin{equation}
A\underline{u}^0=\mu\textrm{ and }A\overline{u}^0=\nu\textrm{ in }\Omega.
\end{equation}
Assuming $f_i\in[L^\infty(\Omega)]^N$, we have $\underline{u}^0$, $\overline{u}^0\in C^{1,\alpha}(\overline{\Omega})$ (see [20]) and we may define the constant $\lambda_0$ by
$$
\lambda_0=\max_{\overline{\Omega}}\overline{u}^0-\min_{\overline{\Omega}}\underline{u}^0>0.
$$
\begin{teo}
Let $f_i\in[L^\infty(\Omega)]^N$ and suppose 
$$\varphi_{ij}+\psi_{ik}\geq\lambda_0>0,\textrm{ }\forall i,j,k=1,\ldots,N.
$$
Then $(5.5)$ admits at least a maximal solution $\overline{\textbf{u}}$ and a minimal solution $\underline{\textbf{u}}$.
\end{teo}

\textbf{Proof.} Let $\mu$, $\nu\in L^\infty(\Omega)$ be given by $(5.6)$. By regularity, the set
$$
\mathbb{D}:=\{\textbf{v}=(v_1,\ldots,v_N)\in [W_0^{1,G}(\Omega)]^N\textrm{ }\mu\leq Av_i\leq\nu\textrm{ in }\Omega,\textrm{ }i=1,\ldots,N\}
$$
is bounded in $[C^{1,\alpha}(\overline{\Omega})]^N$ for some $\alpha>0$. By comparison we have $\underline{u}^0\leq v_i\leq\overline{u}^0$ in $\Omega$, $i=1,\ldots,N$, for every $\textbf{v}\in\mathbb{D}$ and, of course, $\underline{\textbf{u}}^0=(\underline{u}^0,\ldots,\underline{u}^0)$ and $\overline{\textbf{u}}^0=(\overline{u}^0,\ldots,\overline{u}^0$ also belong to $\mathbb{D}$. Therefore, if $\textbf{v}\in\mathbb{D}$ we have
$$
\max_{\overline{\Omega}}v_j-\min_{\overline{\Omega}}v_k\leq\lambda_0,\textrm{ }\forall j,k
$$
and
$$
\mathbb{K}(\textbf{v})=\{\textbf{z}=(z_1,\ldots,z_N)\in[W_0^{1,G}(\Omega)]^N:\textrm{ }\Psi_i(\textbf{v})\leq z_i\leq\Phi_i(\textbf{v}),\textrm{ }i=1,\ldots,N\}\neq\emptyset.
$$
If we denote by $\textbf{w}=\sigma(\textbf{v})$ the unique solution of
\begin{equation}
\textbf{w}\in\mathbb{K}(\textbf{v}):\textrm{ }\langle A\textbf{w}-\textbf{L},\textbf{z}-\textbf{w}\rangle\geq0,\textrm{ }\forall \textbf{z}\in\mathbb{K}(\textbf{v}),
\end{equation}
by Lewy-Stampacchia inequalities and Theorem 4.1, it satisfies
$$
\bigg(\bigwedge_{i\neq j}Av_j\bigg)\wedge f_i\leq Aw_i\leq\bigg(\bigvee_{i\neq j}Av_j\bigg)\vee f_i\textrm{ in }\Omega.
$$
Consequently, $\textbf{w}\in\mathbb{D}$ and $\sigma(\mathbb{D})\subset\mathbb{D}$.

Using $(5.8)$ we define $\overline{\textbf{u}}^1=\sigma(\overline{\textbf{u}}^0)$ and by the weak maximum principle, we observe that $\overline{u}_i^1\leq\overline{u}^0$ in $\Omega$, for all $i=1,\ldots,N$. By iteration, let $\overline{\textbf{u}}^{m+1}=\sigma(\overline{\textbf{u}}^m)$ for $m=1,2,\ldots$ and observe that, if componentwise $\overline{u}_i^m\leq\overline{u}_i^{m-1}$ then
$$
\Phi_i(\overline{\textbf{u}}^m)\leq\Phi_i(\overline{\textbf{u}}^{m-1})\textrm{ and }
\Psi_i(\overline{\textbf{u}}^m)\leq\Psi_i(\overline{\textbf{u}}^{m-1}).
$$
By monotonicity with respect to the obstacles, we conclude then
$$
\overline{u}_i^{m+1}\leq\overline{u}_i^m,\textrm{ }i=1,\ldots,N\textrm{ for all }m=0,1,2,\ldots
$$
Therefore we have constructed a decreasing (and bounded in $[C^{1,\alpha}(\overline{\Omega})]^N$) sequence $\overline{\textbf{u}}^m$ in $\mathbb{D}$, which converges as $m\rightarrow\infty$ to $\overline{\textbf{u}}\in\mathbb{D}$ in $[C^1(\overline{\Omega})]^N$. In particular, this implies, $\overline{\textbf{u}}\in\mathbb{K}(\overline{\textbf{u}})$ and
$$
\Phi(\overline{\textbf{u}}^m)\rightarrow\Phi(\overline{\textbf{u}})\textrm{ and }
\Psi(\overline{\textbf{u}}^m)\rightarrow\Psi(\overline{\textbf{u}})\textrm{ in }[W^{1,G}(\Omega)]^N.
$$
Since $\overline{\textbf{u}}^{m+1}\in\mathbb{K}(\overline{\textbf{u}}^m)$ also
$$
\langle A\overline{\textbf{u}}^{m+1}-\textbf{L},\overline{\textbf{v}}-\overline{\textbf{u}}^{m+1}\rangle\geq0,\textrm{ }\forall\overline{\textbf{v}}\in\mathbb{K}(\overline{\textbf{u}}^m).
$$
By Theorem 3.6, $\overline{\textbf{u}}$ solves $(5.5)$.

Analogously, for $m=0,1,2,\ldots$ we construct an increasing sequence
$$
\underline{\textbf{u}}^m\leq\underline{\textbf{u}}^{m+1}=\sigma(\underline{\textbf{u}}^m)\rightarrow\underline{\textbf{u}}\textrm{ in }[C^1(\overline{\Omega})]^N\cap\mathbb{D}
$$
and easily conclude that $\underline{\textbf{u}}\in\mathbb{K}(\underline{\textbf{u}})$ also solves $(5.5)$.

If $\textbf{u}=(u_1,\ldots,u_N)$ is any other solution to $(5.5)$, by the Lewy-Stampacchia inequalities we have for all $i=1,\ldots,N$
$$
\mu\leq\bigg(\bigwedge_{i\neq j}Au_j\bigg)\wedge f_i\leq Au_i\leq\bigg(\bigvee_{i\neq j}Au_j\bigg)\vee f_i\leq\nu
$$
and so, by comparison, $\underline{u}^0\leq u_i\leq\overline{u}^0$. Hence $\underline{u}_i^1\leq u_i\leq\overline{u}_i^1$, by monotonicity of $\Psi$ and $\Phi$, and so also $\underline{u}_i^m\leq u_i\leq\overline{u}_i^m$ for all $m$ by recurrence.

We then conclude $\underline{u}_i\leq u_i\leq\overline{u}_i$, $i=1,\dots,N$ that yields the minimality of $\underline{\textbf{u}}$ and maximality of $\overline{\textbf{u}}$.$\blacksquare$
\begin{rem}
The uniqueness of the solution to $(5.5)$ is an open problem.
\end{rem}

\end{document}